\documentclass[a4paper,10pt]{amsart}
\usepackage{amssymb}
\usepackage{amsmath}
\usepackage{amsthm}
\usepackage[pdftex]{graphicx}
\usepackage[english]{babel}
\usepackage{array}
\usepackage{titletoc}
\usepackage[T1]{fontenc}
\usepackage[utf8]{inputenc}
\usepackage[all]{xy}
\usepackage[vcentermath]{youngtab}
\usepackage{ytableau}
\usepackage{stmaryrd}

\usepackage[mathcal]{euscript}

\usepackage[bbgreekl]{mathbbol}

\usepackage{enumerate}

\usepackage{hyperref}

\allowdisplaybreaks[1]

 \numberwithin{equation}{section}

  \theoremstyle{plain}
\newtheorem{theoreme}[equation]{Theorem}
\newtheorem{lemme}[equation]{Lemma}
\newtheorem{lemme-def}[equation]{Lemma-Definition}
\newtheorem{proposition}[equation]{Proposition}
\newtheorem{corollaire}[equation]{Corollary}

  \theoremstyle{definition}
\newtheorem{definition}[equation]{Definition}
\newtheorem{conj}[equation]{Conjecture}

  \theoremstyle{remark}
\newtheorem{remarque}[equation]{Remark}
\newtheorem{notations}[equation]{Notations}
\newtheorem{exemple}[equation]{Example}

 \newcommand{\theo}{\begin{theoreme}}
 \newcommand{\defi}{\begin{definition}}
 \newcommand{\rema}{\begin{remarque}}
 \newcommand{\prop}{\begin{proposition}}
 \newcommand{\coro}{\begin{corollaire}}
 \newcommand{\lemm}{\begin{lemme}}
 \newcommand{\exem}{\begin{exemple}}
 \newcommand{\nota}{\begin{notations}}

 \newcommand{\etheo}{\end{theoreme}}
 \newcommand{\edefi}{\end{definition}}
 \newcommand{\erema}{\end{remarque}}
 \newcommand{\eprop}{\end{proposition}}
 \newcommand{\ecoro}{\end{corollaire}}
 \newcommand{\elemm}{\end{lemme}}
 \newcommand{\eexem}{\end{exemple}}
 \newcommand{\enota}{\end{notations}}

 \newcommand{\becen}{\begin{center}}
 \newcommand{\ecen}{\end{center}}

 \newcommand{\benu}{\begin{enumerate}}
 \newcommand{\eenu}{\end{enumerate}}

 \newcommand{\bite}{\begin{itemize}}
 \newcommand{\eite}{\end{itemize}}

\def\[#1\]{\begin{align*}#1\end{align*}}

\newcommand{\demo}{\begin{proof}}
\newcommand{\edemo}{\end{proof}}

\newcommand{\RR}{\mathbb R}
\newcommand{\CC}{\mathbb C}
\newcommand{\QQ}{\mathbb Q}
\newcommand{\ZZ}{\mathbb Z}
\newcommand{\NN}{\mathbb N}

\newcommand{\Hom}{\operatorname{Hom}}

\newcommand{\aaa}{\alpha}
\newcommand{\bbb}{\beta}

\newcommand{\Irr}{\operatorname{Irr}}

\newcommand{\rarrow}{\rightarrow}

\newcommand{\Rarrow}{\Rightarrow}

\newcommand{\isom}{\overset{\sim}{\rightarrow}}

\newcommand{\ima}{\operatorname{Im}}

\usepackage{lipsum}

\author{Tristan Bozec\\ \\    with an appendix by Anton Mellit}

\title[The global nilpotent cone]{Irreducible components of the global nilpotent cone}
\date{}

\begin{document}

\maketitle
\begin{abstract}
This paper gives a combinatorial description of the set of irreducible components of the semistable locus of the global nilpotent cone, in genus $\ge2$.

\noindent{\emph{Keywords:} Moduli stacks of sheaves over curves of higher genus, Higgs sheaves, Lagrangian subvarieties, stability.

}
\end{abstract}
\tableofcontents

\section*{Introduction}

Given a smooth projective curve $X$ of genus $g$, the moduli stack of Higgs sheaves of rank $r$ and degree $d$ is known to be of dimension $2(g-1)r^2$. It can be viewed as the cotangent stack of the stack of coherent sheaves of class $(r,d)$ over $X$, and Laumon proved in~\cite{Laumon} that the substack $\mathbf\Lambda_{r,d}$ of nilpotent Higgs pairs is Lagrangian (see also~\cite{Faltings,GinzLag}). This substack, which is the $0$-fiber of the Hitchin map, is a global analog of the nilpotent cone and a plays a critical role in the geometric Langlands program. The global nilpotent cone is highly singular, and one first interesting step toward its comprehension is the study of its set of irreducible components (see~\cite[2.10.3]{BeilinsonDrinfeld} for rather implicit results in this direction).

The stack of \textit{stable} Higgs pairs is known to be smooth, and several results have been proved recently regarding the counting of the number of stable irreducible components of the global nilpotent cone. Its is known from~\cite[Corollary 3.11]{MR2166085} that the Poincar\'e polynomial of genus $g$ twisted character varieties, and hence of the diffeomorphic moduli spaces of stable Higgs bundles, is independent from the degree
$d$, provided that it is coprime to the rank $r$. In~\cite{MR2453601}, Hausel and Rodriguez-Villegas establish several conjectures dealing with the $E$-polynomial (a specialization of the mixed Hodge polynomial) of these character varieties. In particular, they conjecture a combinatorial relation with the Kac polynomial $A_{g,r}$ of the quiver with one vertex and $g$ loops in dimension $r$, which counts absolutely indecomposable isoclasses of $g$-tuples of matrices over finite fields.

With a different perspective, Schiffmann establishes in~\cite{SchiffAnn} that the number of absolutely indecomposable vector bundles of rank $r$ and degree $d$ over $\mathbb F_q$ (still over a curve $X$ of genus $g$) is given by an expression $A_{g,r,d}$, polynomial in the Weil numbers of $X$. These polynomials are therein proved to be related to the moduli space of stable Higgs bundles, and, for instance, the number of stable irreducible components of $\mathbf\Lambda_{r,d}$ is given by $A_{g,r,d}(0)$.

In a recent work~\cite{Mellit}, Mellit relates the formulas obtained in~\cite{MR2453601,SchiffAnn}, and proves as a consequence that the polynomials $A_{g,r,d}$ and the $E$-polynomials of the moduli spaces of stable Higgs pairs are both independent of the degree $d$ coprime with $r$. A very particular consequence of this work is the equality $A_{g,r,d}(0)=A_{g,r}(1)$.

The aim of the present paper is to give a combinatorial description of the set of irreducible components of $\mathbf\Lambda_{r,d}$, and explain which ones subsist in the subset of semistable components. It is motivated by the $W=P$ conjecture claimed by de Cataldo, Hausel and Migliorini~\cite{MR2912707}. In the light of the works above-mentioned, one can expect the polynomial $A_{g,r}$ to play a role in the understanding of the perverse filtration. Adding structure to the set of irreducible components could lead to an interpretation of each of the coefficients of $A_{g,r}$, rather than just their sum ($=A_{g,r}(1)$).

The first main result of this paper is Corollary~\ref{result1}, which states that the set of irreducible components of the global nilpotent cone is given by the very natural decomposition in twisted Jordan cells, which are smooth. It is based on a direct computation of the dimensions of these cells and previous works~\cite{SchiffAnn,MS}. Then we move on to the semistable locus and obtain Theorem~\ref{result2} which gives purely combinatorial conditions on the twisted Jordan type to be semistable. The proof uses an analogous result from~\cite{BGGH} obtained in the context of moduli stacks of chains, and shows that semistability can be tested on the most `simple' subsheaves - the ones built with iterated kernels and images. The proof is constructive and do not rely on the coprimality of $r$ and $d$, in particular we get in Corollary~\ref{corocell} that the attracting cells are irreducible in any case.

The Corollary~\ref{polytop} describes this set of semistable irreducible components in terms of integral polytopes, which sheds a new light on the quantity $A_{g,r}(1)$, whose behaviour is still very poorly understood.

\bigskip

\noindent{\textit{Acknowledgements.}}
This work was started during the author's postdoctoral appointment at MIT, before being supported by the LABEX MILYON (ANR-10-LABX-0070) of Universit\'e de Lyon, within the program ``Investissements d'Avenir'' (ANR-11-IDEX-0007) operated by the French National Research Agency (ANR). The author would like to thank Olivier Schiffmann for numerous extremely useful discussions.

\section{Recollection on coherent and Higgs sheaves}

\subsection{Coherent sheaves over a curve}
Let $X$ be a smooth projective curve of genus $g$ over a field $\mathbb k$.
We will denote by $\text{\sffamily Coh}$ the category of coherent sheaves over $X$, and by\[
[\mathcal F]=(\operatorname{rank}\mathcal F,\deg\mathcal F)\in\mathbb H=\{(r,d)\in\NN\times\ZZ\mid d\ge0\text{ if }r=0\}\]
the class of $\mathcal F\in\text{\sffamily Coh}$. We will denote by $\text{\sffamily Coh}_{r,d}\subset\text{\sffamily Coh}$ the subcategory of coherent sheaves of class $(r,d)$. If $\aaa=(r,d)$, we write $r=\operatorname{rank}\aaa$ and $d=\deg\aaa$. For any $\aaa=(r,d)\in\mathbb H$ and $p\in\ZZ$, we set $\aaa(p)=(r,d+pr)$ so that if $\mathcal F\in\text{\sffamily Coh}$ and $D$ is a divisor of degree $p$ over $X$ we have \[
[\mathcal F(D)]=[\mathcal F\otimes D]=[\mathcal F](p).\]
We will use the usual \emph{slope} defined on $\mathbb H$ by $\mu(r,d)=d/r\in\QQ\cup\{\infty\}$ and we set $\mu(\mathcal F)=\mu([\mathcal F])$. We say that $\mathcal F$ is \emph{semistable} if\[
\{0\}\subset\mathcal G\subset\mathcal F\Rightarrow\mu(\mathcal G)\le\mu(\mathcal F),\]
\emph{stable} if the right-hand side inequality is strict. Note that these notions coincide if $\deg\mathcal F$ and $\operatorname{rank}\mathcal F$ are coprime. We will use the following basic property.

\prop\label{basique}
For any short exact sequence\[
0\rarrow\mathcal E\rarrow\mathcal F\rarrow\mathcal G\rarrow0\]
 in $\textup{\sffamily Coh}$, one of the following is true\[
\mu(\mathcal E)&<\mu(\mathcal F)<\mu(\mathcal G)\\
\mu(\mathcal E)&=\mu(\mathcal F)=\mu(\mathcal G)\\
\mu(\mathcal E)&>\mu(\mathcal F)>\mu(\mathcal G).\]
\eprop

The category $\textup{\sffamily Coh}$ being hereditary (\textit{i.e.}\ of homological dimension one or less), the Euler form is defined by\[
\langle\mathcal F,\mathcal G\rangle=\dim\Hom(\mathcal F,\mathcal G)-\dim\operatorname{Ext}^1(\mathcal F,\mathcal G),\]
and we will denote by $(-,-)$ its symmetrized version\[
(\mathcal F,\mathcal G)=\langle \mathcal F,\mathcal G\rangle+\langle\mathcal G,\mathcal F\rangle.\]  The Euler form only depends on the class of the sheaves, and satisfies\[
\langle(r,d),(r',d')\rangle=(1-g)rr'+rd'-r'd\]
thanks to the Riemann--Roch theorem.

\subsection{Higgs sheaves}

A \emph{Higgs sheaf} is a pair $(\mathcal F,\theta)$, where $\mathcal F\in\textup{\sffamily Coh}$ and $\theta\in \Hom(\mathcal F,\mathcal F(\Omega))$, $\Omega$ being the canonical divisor of degree $l=2g-2$. We will denote by ${\mathbf M}_{r,d}$ the moduli stack of pairs $(\mathcal F,\theta)$ satisfying $[\mathcal F]=(r,d)$, whose dimension is $lr^2$. A Higgs sheaf $(\mathcal F,\theta)$ is said to be semistable\[
\left.\begin{aligned}\{0\}\subset\mathcal G\subset\mathcal F~\\
\theta(\mathcal G)\subseteq\mathcal G(\Omega)~\end{aligned}\right\}
\Rightarrow\mu(\mathcal G)\le\mu(\mathcal F),\]
stable if the right-hand side inequality is strict.
Semistability defines an open substack $\mathbf M^\textup{sst}_{r,d}\subset \mathbf M_{r,d}$.

\subsection{The global nilpotent cone}

\defi For any $(\mathcal F,\theta)\in{\mathbf M}_{r,d}$ and $k\ge1$, set \[
\theta^k=\theta((k-1)\Omega)\circ\dots\circ\theta(\Omega)\circ\theta:\mathcal F\rarrow\mathcal F(k\Omega).\]
 A pair $(\mathcal F,\theta)$ is said to be \emph{nilpotent} if $\theta^k=0$ for some $k$, and we denote by\[
\boldsymbol\Lambda_{{r,d}}=\{(\mathcal F,\theta)\in{\mathbf M}_{r,d}\mid(\mathcal F,\theta)\text{ nilpotent}\}\]
the \emph{global nilpotent cone}.
\edefi

It is nothing but the zero fiber of the Hitchin map $\mathbf M_{r,d}\rarrow\oplus_{1\le i \le r}H^0(X,\Omega^i)$, mapping $(\mathcal F,\theta)$ to the coefficients of the characteristic polynomial of $\theta$. It is known, thanks to Laumon~\cite{Laumon}, to be a Lagrangian substack of $\mathbf M_{r,d}$, but its irreducible components are still not well understood (see [Ibid., Remarque 3.9 (ii)]). The aim of this article is to give a precise combinatorial description of these components, as well as the ones of the semistable locus $\boldsymbol\Lambda_{r,d}^\textup{sst}$.

\section{The twisted Jordan cell decomposition}

\subsection{The setting}\label{setting}
In this section we will recall and make use of notations and results established in~\cite{MS,SchiffAnn}.
Consider a nilpotent Higgs sheaf $(\mathcal F,\theta)\in\boldsymbol\Lambda_{r,d}$, with $r>0$ but $r$ and $d$ not necessarily coprime.
Set $\mathcal F_k=(\ima\theta^k)(-k\Omega)$ and denote by $s$ the nilpotency index of $\theta$.
We have a chain of epimorphisms\[
\mathcal F_0\twoheadrightarrow\mathcal F_1(\Omega)\twoheadrightarrow\dots\twoheadrightarrow\mathcal F_s(s\Omega)=\{0\}\]
that allows us to define $\mathcal F'_k=\ker\{\mathcal F_k\rarrow\mathcal F_{k+1}(\Omega)\}$.
We also have a chain of inclusions\[
\{0\}=\mathcal F_s\subset\mathcal F_{s-1}\subset\dots\subset\mathcal F_1\subset\mathcal F_0=\mathcal F\]
whose successive quotients are denoted by $\mathcal F''_k=\mathcal F_k/\mathcal F_{k+1}$. These two chains induce the following ones\[
\mathcal F''_0\twoheadrightarrow\mathcal F''_1(\Omega)\twoheadrightarrow\dots\twoheadrightarrow\mathcal F''_s(s\Omega)=\{0\}\]\[
\{0\}=\mathcal F'_s\subset\mathcal F'_{s-1}\subset\dots\subset\mathcal F'_1\subset\mathcal F'_0,\]
and we define\[
\aaa_k=(r_k,d_k)=\left[\ker\{\mathcal F''_{k-1}((k-1)\Omega)\rarrow\mathcal F''_k(k\Omega)\}\right].\]
The family $\boldsymbol\aaa=(\mathbf r,\mathbf d)=(r_1,\dots,r_s,d_1,\dots,d_s)$ is called the \emph{Jordan type} of $(\mathcal F,\theta)$ and is denoted by $J(\mathcal F,\theta)$. We will call $s$ the \emph{length} of $\boldsymbol\aaa$ or $\mathbf r$. One good way to understand its definition is to fill the triangular Young tableau $T_s$ of size $s$ in the following way (here with s=4) \begin{align}\label{TS}
\ytableausetup
 {mathmode, boxsize=2.8em,centertableaux}
\begin{ytableau}
\scriptstyle{\aaa_4}   \\
\scriptstyle\aaa_4(-l) & \scriptstyle\aaa_3 \\
\scriptstyle{\aaa_4(-2l)} &\scriptstyle \aaa_3(-l) & \scriptstyle\aaa_2\\
\scriptstyle{\aaa_4(-3l)} &\scriptstyle \aaa_3(-2l) & \scriptstyle\aaa_2(-l)&\scriptstyle\aaa_1
  \end{ytableau}
\end{align}
and then notice that \[
[\mathcal F''_k]=\sum_{i>k}\aaa_i(-kl)\] is the sum of the classes in the boxes of the $k$-th subdiagonal. Hence, $[\mathcal F_k]=\sum_{i\ge k}[\mathcal F''_i]$ corresponds to the region below this subdiagonal. Denote by $[R]$ the sum of the classes of the boxes in a region $R$. We have, for instance with $s=5$,\[
[\mathcal F''_1]=\left[~\ytableausetup
 {smalltableaux,centertableaux}
\begin{ytableau}
  \\
*(black) &  \\
& *(black)  & \\
&& *(black)  & \\
&&& *(black)  &
  \end{ytableau}~ \right]
\qquad
[\mathcal F_2]=\left[~\ytableausetup
 {smalltableaux,centertableaux}
\begin{ytableau}
  \\
 &  \\
*(black)&   & \\
*(black)& *(black)&   & \\
*(black)&*(black)&*(black)&   &
  \end{ytableau}~\right]
\qquad
[\mathcal F''_4]=[\mathcal F_4]=\left[~\ytableausetup
 {smalltableaux,centertableaux}
\begin{ytableau}
  \\
 &  \\
&   & \\
 &&&\\
 *(black)&&&&
  \end{ytableau}~\right]
\]
where the classes are summed over the blackened regions.
In particular the sum over all boxes is $(r,d):=[\mathcal F]$ and we write $\boldsymbol\aaa\vdash(r,d)$ (and call $(r,d)$ the \textit{size} of $\boldsymbol\aaa$). This implies the following equality\begin{align}\label{deun}
\sum_kkd_k=d+l\sum_k\frac{k(k-1)}{2}r_k.\end{align}
We also have \[
[\ker\theta^{k}]-[\ker\theta^{k-1}]&=[\mathcal F_0]-[\mathcal F_k(k\Omega)]-\left\{[\mathcal F_0]-[\mathcal F_{k-1}((k-1)\Omega)]\right\}\\
&=[\mathcal F_{k-1}((k-1)\Omega)]-[\mathcal F_k(k\Omega)]\\
&=\sum_{j>i\ge k-1}\aaa_j((k-1-i)l)-\sum_{j>i\ge k}\aaa_j((k-i)l)\\
&=-\sum_{j>k-1}\aaa_j(l(j-k+1)(j-k)/2)\\&\qquad\qquad\qquad+\sum_{j>k}\aaa_j(l(j-k)(j-k-1)/2)\\
&=\sum_{j\ge k}\aaa_j((k-j)l)\]
which corresponds to the $k$-th (from the bottom) horizontal strip. Graphically, we have for instance\[
[\ker\theta]=\left[~\ytableausetup
 {smalltableaux,centertableaux}
\begin{ytableau}
  \\
 &  \\
&  & \\
&&   & \\
*(black)&*(black)&*(black)& *(black)  & *(black)
  \end{ytableau}~\right]
\qquad
[\ker\theta^3]=\left[~\ytableausetup
 {smalltableaux,centertableaux}
\begin{ytableau}
  \\
 &  \\
*(black)& *(black)  & *(black)\\
*(black)& *(black)&  *(black) & *(black)\\
*(black)&*(black)&*(black)&  *(black) & *(black)
  \end{ytableau}~\right]
\]
We call \emph{canonical} the subsheaves of $\mathcal F$ obtained by intersections and sums of the $\mathcal F_k$ and the $\ker\theta^k$. The corresponding regions in the Young tableau are the ones saturated in the west, south and south-east directions (note that there is a mistake in the corresponding statement in~\cite[3.1]{SchiffAnn}). We denote by $\mathcal R$ this set of regions, which we will also call canonical. The slope of the sheaf $\mathcal F_R$ corresponding to a region $R\in\mathcal R$ is given by\[
\mu^{\boldsymbol\aaa}(R)=\dfrac{ d_R}{ r_R}=\dfrac{\sum_{\blacksquare\in R}\deg\blacksquare}{\sum_{\blacksquare\in R}\operatorname{rank}\blacksquare}
\]
with respect to the filling~(\ref{TS}) by $\boldsymbol\aaa=(\mathbf r,\mathbf d)=J(\mathcal F,\theta)$ (we will call $\mu^{\boldsymbol\aaa}(R)$ the \emph{$\boldsymbol\aaa$-slope} of $R$).
For instance\[
\mu(\ker\theta^2\cap\mathcal F_2+\ker\theta\cap\mathcal F_1)&=\mu^{\boldsymbol\aaa}\!\!\left(\ytableausetup
 {smalltableaux,centertableaux}
\begin{ytableau}
  \\
 &  \\
&&\\
*(black)& *(black)&  &\\
*(black)&*(black)&*(black)&  *(black) &
  \end{ytableau}\right)\\
&=\dfrac{2d_5+2d_4+d_3+d_2-7lr_5-5lr_4-2lr_3-lr_2}{2r_5+2r_4+r_3+r_2}.
\]
The definition of the Jordan type yields a cell decomposition\[
\boldsymbol\Lambda_{r,d}=\bigsqcup_{\boldsymbol\aaa\vdash(r,d)}\boldsymbol\Lambda_{\boldsymbol\aaa}\]
where $\boldsymbol\Lambda_{\boldsymbol\aaa}=J^{-1}(\boldsymbol\aaa)\subseteq\boldsymbol\Lambda_{r,d}$.

\begin{exemple}
For $s=4$, the set $\mathcal R$ of canonical regions consists in\[
\ytableausetup
 {smalltableaux,centertableaux}
\begin{ytableau}
   \\
&\\
&&  \\
*(black)&&&
  \end{ytableau};
\begin{ytableau}
   \\
&\\
&&  \\
*(black)&*(black)&&
  \end{ytableau};
  \begin{ytableau}
   \\
&\\
&&  \\
*(black)&*(black)&*(black)&
  \end{ytableau};
  \begin{ytableau}
   \\
&\\
&&  \\
*(black)&*(black)&*(black)&  *(black)
  \end{ytableau};
  \begin{ytableau}
   \\
&\\
*(black)&&  \\
*(black)&*(black)&&
  \end{ytableau};
  \begin{ytableau}
   \\
&\\
*(black)&&  \\
*(black)&*(black)&*(black)&
  \end{ytableau};
  \begin{ytableau}
   \\
&\\
*(black)&&  \\
*(black)&*(black)&*(black)&  *(black)
  \end{ytableau};
    \begin{ytableau}
   \\
&\\
*(black)& *(black)&  \\
*(black)&*(black)&*(black)&
  \end{ytableau};
\]
  \[
\ytableausetup
 {smalltableaux,centertableaux}
  \begin{ytableau}
   \\
&\\
*(black)& *(black)&  \\
*(black)&*(black)&*(black)&  *(black)
  \end{ytableau};
    \begin{ytableau}
   \\
&\\
*(black)& *(black)&*(black)  \\
*(black)&*(black)&*(black)&  *(black)
  \end{ytableau}
 ;
     \begin{ytableau}
   \\
*(black) &\\
*(black)& *(black)& \\
*(black)&*(black)&*(black)&
  \end{ytableau};
    \begin{ytableau}
   \\
*(black) &\\
*(black)& *(black)& \\
*(black)&*(black)&*(black)&  *(black)
  \end{ytableau};
   \begin{ytableau}
   \\
*(black) &\\
*(black)& *(black)& *(black)  \\
*(black)&*(black)&*(black)&  *(black)
  \end{ytableau};
   \begin{ytableau}
   \\
*(black) &*(black)\\
*(black)& *(black)& *(black)  \\
*(black)&*(black)&*(black)&  *(black)
  \end{ytableau};
   \begin{ytableau}
*(black)   \\
*(black) &*(black)\\
*(black)& *(black)& *(black)  \\
*(black)&*(black)&*(black)&  *(black)
  \end{ytableau}.\]
\end{exemple}

\subsection{Irreducible components}

In this section, we will study the map of stacks \[
\left.\begin{aligned}\pi_{\boldsymbol\aaa}\colon\boldsymbol\Lambda_{\boldsymbol\aaa}&\longrightarrow&&\prod_k\textbf {Coh}_{\aaa_k}\\
(\mathcal F,\theta)&\longmapsto&&\big(\ker\{\mathcal F''_{k-1}((k-1)\Omega)\rarrow\mathcal F''_k(k\Omega)\}\big)_{\!k}\end{aligned}\right.\]
for any Jordan type $\boldsymbol\aaa$. Denote by $\mathbf F_{\boldsymbol\aaa}$ the stack of chains of epimorphisms\[
\mathcal H_0\twoheadrightarrow\mathcal H_1\twoheadrightarrow\dots\twoheadrightarrow\mathcal H_s=\{0\}
\]
satisfying $\aaa_k=\left[\ker\{\mathcal H_{k-1}\twoheadrightarrow\mathcal H_k\}\right]$
and write $\pi_{\boldsymbol\aaa}=\rho_{\boldsymbol\aaa}\circ\chi_{\boldsymbol\aaa}$ where \[
\left.\begin{aligned}\chi_{\boldsymbol\aaa}\colon\boldsymbol\Lambda_{\boldsymbol\aaa}&\longrightarrow&&\mathbf F_{\boldsymbol\aaa}\\
(\mathcal F,\theta)&\longmapsto&&\left(\mathcal F''_0\twoheadrightarrow\mathcal F''_1(\Omega)\twoheadrightarrow\dots\twoheadrightarrow\mathcal F''_s(s\Omega)\right)\end{aligned}\right.\]
and\[
\left.\begin{aligned}\rho_{\boldsymbol\aaa}\colon\mathbf F_{\boldsymbol\aaa}&\longrightarrow&&\prod_k\textbf {Coh}_{\aaa_k}\\\left(
\mathcal H_0\twoheadrightarrow\mathcal H_1\twoheadrightarrow\dots\twoheadrightarrow\mathcal H_s\right)&\longmapsto&&\big(\ker\{\mathcal H_{k-1}\rightarrow\mathcal H_k\}\big)_{\!k}.\end{aligned}\right.\]

The following equalities are obtained in~\cite[Proposition 5.2]{MS} and ~\cite[3.1]{SchiffAnn} respectively.
\prop\label{recollschiff} The maps $\chi_{\boldsymbol\aaa}$ and $\rho_{\boldsymbol\aaa}$ are iterations of vector bundle stacks and their respective relative dimensions are\[
d_{\chi_{\boldsymbol\aaa}}&=-\sum_k\langle \mathcal F''_k,\mathcal F'_{k+1}\rangle\\
d_{\rho_{\boldsymbol\aaa}}&=-\sum_{i<j}\langle\aaa_j,\aaa_i\rangle.\]
\eprop

\coro\label{result1} The set of irreducible components of $\boldsymbol\Lambda_{r,d}$ is \[
\Irr\boldsymbol\Lambda_{r,d}=\left\{\overline{\boldsymbol\Lambda_{\boldsymbol\aaa}}\mid\boldsymbol\aaa\vdash(r,d)\right\}.\]\ecoro

\demo From Proposition~\ref{recollschiff}, since all $\textbf{Coh}_{\alpha_k}$ are irreducible, we know that the $\boldsymbol\Lambda_{\boldsymbol\aaa}$ are also irreducible (even smooth), thus we just have to prove that they have the same dimension. From Proposition~\ref{recollschiff} and the fact that $\dim\textbf{Coh}_{r,d}=-\langle(r,d),(r,d)\rangle$, we get\[
\dim\boldsymbol\Lambda_{\boldsymbol\aaa}&=-\sum_k\langle \mathcal F''_k,\mathcal F'_{k+1}\rangle-\sum_{i<j}\langle\aaa_j,\aaa_i\rangle-\sum_k\langle\aaa_k,\aaa_k\rangle\\&=-\sum_{i>k}\langle \aaa_i(-kl),\mathcal F'_{k+1}\rangle-\sum_{i\le j}\langle\aaa_j,\aaa_i\rangle.
\]
Now, since $[\ker\theta^k]=\sum_{j\le k-1}[\mathcal F'_j(j\Omega)]$, we have\[
[\mathcal F'_k]=([\ker\theta^{k+1}]-[\ker\theta^k])(-kl)=\sum_{j>k}\aaa_j((1-j)l),
\]
thus\[
\dim\boldsymbol\Lambda_{\boldsymbol\aaa}&=-\sum_{\substack{i>k\\j>k+1}}\langle\aaa_i(-kl),\aaa_j((1-j)l)\rangle-\sum_{i\le j}\langle\aaa_j,\aaa_i\rangle\\
&=-\sum_{\substack{i>k\\j>k+1}}(\langle\aaa_i,\aaa_j\rangle+l(k+1-j)r_ir_j)-\sum_{j\le i}\langle\aaa_i,\aaa_j\rangle\\
&=-\sum_{i\ge j}\left(j\langle\aaa_i,\aaa_j\rangle-l\frac{j(j-1)}{2}r_ir_j\right)\\
&\qquad\qquad-\sum_{i< j}\left(i\langle\aaa_i,\aaa_j\rangle+l\frac{i(i+1-2j)}{2}r_ir_j\right)\\
&=-\sum_{i<j}i(\aaa_i,\aaa_j)-\sum_i\left(i\langle \aaa_i,\aaa_i\rangle-l\frac{i(i-1)}{2}r_i^2\right)\\
&\qquad\qquad+\frac{l}{2}\sum_{i<j}(i(i-1)-i(i+1-2j))r_ir_j\\
&=l\sum_{i<j}ir_ir_j+\frac{l}{2}\sum_i(i+i(i-1))r_i^2+l\sum_{i<j}i(j-1)r_ir_j\\
&=(g-1)r^2
\]
as expected.
\edemo

\section{The semistable locus}

As semistability is an open condition, we know that the irreducible components of $\boldsymbol\Lambda_{(r,d)}^\textup{sst}$ form a subset of $\Irr\boldsymbol\Lambda_{(r,d)}$ which can be now identified, thanks to Corollary~\ref{result1}, with the set of all Jordan types of size ${(r,d)}$. We say that a Jordan type is semistable if it appears in $\Irr\boldsymbol\Lambda_{{(r,d)}}^\textup{sst}$. The aim of this section is to prove the following.

\theo\label{result2} Assume that $g\ge2$. A Jordan type $\boldsymbol\aaa\vdash(r,d)$ is semistable if and only if for every canonical strict subregion $R\in\mathcal R$ we have\begin{align}\label{slope}
\mu^{\boldsymbol\aaa}(R)\le\frac{d}{r}.\end{align}\etheo

 Note that this condition is strictly numerical and obviously necessary since a canonical subsheaf $\mathcal G$ satisfies by construction $\theta(\mathcal G)\subseteq\mathcal G(\Omega)$. This result is optimal in the way that it says that it is sufficient to test (generic) semistability on the most trivial $\theta$-stable subsheaves. To prove it, we will use the results of~\cite{BGGH}, which deals with the moduli stack of chains $ \mathcal E_{s}\rarrow\dots\rarrow\mathcal E_1$. In this article is obtained an analogous result in the way that it gives necessary and sufficient conditions on the numerical invariants $(n_\bullet,p_\bullet)=[\mathcal E_\bullet]$ of chains for these to be generically semistable (we will call semistable types of chains such invariants).
 The conditions obtained therein are  somewhat unnatural, in the way that they do not correspond to proper subchains - see~[Ibid., Remark 2.11]. However, in this section (\textit{c.f.}\ Proposition~\ref{1chain}), we will build a injection between Jordan types satisfying~(\ref{slope}) and semistable types of chains, which are related to the semistable components of the nilpotent cone in the following way. We consider stability of chains with respect to the $\alpha^\text{Higgs}$-slope \[
 \mu_{\alpha^\text{Higgs}}(\mathcal E_\bullet)=\dfrac{\sum_{1\le i\le s}(\operatorname{deg}(\mathcal E_i)+\aaa^\text{Higgs}_i\operatorname{rank}(\mathcal E_i))}{\sum_{1\le i\le s}\operatorname{rank}(\mathcal E_i)}\]
 where $\alpha^\text{Higgs}=((i-1)(2g-2))_{1\le i\le s}$. Then the direct sum $\oplus_k\mathcal E_k((k-1)\Omega)$ of a semistable chain yields a Higgs pair which is a semistable fixed point under the action $t.(\mathcal F,\theta)=(\mathcal F,t\theta)$ of $\CC^*$ on $\mathbf M_{r,d}$, where $\sum_k[\mathcal E_k((k-1)\Omega)]=(r,d)$. Denote by $\mathbf C_{n_\bullet,p_\bullet}$ the substack of fixed points associated to chains of type $(n_\bullet,p_\bullet)$, and by $\mathbf C_{n_\bullet,p_\bullet}^-$ the corresponding attracting variety
 \[\mathbf C_{n_\bullet,p_\bullet}^-=\{(\mathcal F,\theta)\mid\lim_{t\to0}t.(\mathcal F,\theta)\in\mathbf C_{n_\bullet,p_\bullet}\}.\]
 It is known (see e.g.\ \cite[\S6]{BGGH}) that the closures of these attracting varieties, for $(n_\bullet,p_\bullet)$ of semistable type, are unions of irreducible components of $\mathbf\Lambda_{r,d}^\text{sst}$. Hence, building the injection announced will imply Theorem~\ref{result2} as all inequalities in the following chain will have to be equalities:
 \[ \#\Irr\mathbf\Lambda_{r,d}^\text{sst}&\le\#\{\boldsymbol\aaa\vdash(r,d)\mid\text{(3.2)}\}\\
 &\le\#\{\text{semistable types }(n_\bullet,p_\bullet)\mid\textstyle\sum_k(n_k,p_k)((k-1)l)=(r,d)\}\\
 &\le\#\Irr\mathbf\Lambda_{r,d}^\text{sst}.\]
 The last equality can be stated as
 \coro\label{corocell} The closures of the semistable attracting cells $\mathbf C_{n_\bullet,p_\bullet}^-$ are irreducible regardless of the coprimality of $r$ and $d$.\ecoro
 This was only known in the coprime case.

 \bigskip

 Fix for now a type $\boldsymbol\aaa$ of length $s$ satisfying (\ref{slope}), and $(\mathcal F,\theta)\in\boldsymbol\Lambda_{\boldsymbol\aaa}$. We introduce a couple of notions before rephrasing the main result of~\cite{BGGH} in our context (\textit{c.f.}\ Proposition~\ref{recap}).

 \defi\label{1def} We call \textit{$1$-flags} the flags $R_\bullet=(\varnothing=R_0\subset R_1\subset \dots\subset R_t=T_s)$ of subregions of $T_s$ such that:\benu[(i)]
 \item $t=s$;
 \item $R_k\in \mathcal R$;
 \item if the number of boxes in a given column of $T_s$ is increased by $1$ from $R_k$ to $R_{k+1}$, the same must be true for every column on its left.\eenu
We call  \textit{strips} the (noncanonical) subregions $S_k=R_k\setminus R_{k-1}$, and \textit{$1$-chains} the chains $\mathcal E_\bullet=( \mathcal E_{s}\rarrow\dots\rarrow\mathcal E_1)$ associated to a given $1$-flag $R_\bullet$, where $\mathcal E_k=\mathcal F_{R_{k}}/\mathcal F_{R_{k-1}}(-(k-1)\Omega)$, and the morphisms are induced by $\theta$.\edefi

 \rema ~

 \bite
 \item A $1$-chain $\mathcal E_\bullet$ depends on the data $(R_\bullet,\mathcal F,\theta)$, but its type $(n_\bullet,p_\bullet)=[\mathcal E_\bullet]$ only depends on $(R_\bullet,\boldsymbol\aaa)$.
 \item The $1$- notation comes from the fact that $\text{ht}(S_k)=1$ (the height being the number of boxes in the higher column), which is why $\theta$ induces morphisms $\mathcal E_{k+1}\rarrow\mathcal E_k$.
 \item Each strip has a box on the left border of $T_s$.
 \eite
 \erema

\exem ~

 \bite\item The flag of regions ($s=3$)\[
\ytableausetup
 {smalltableaux,centertableaux}
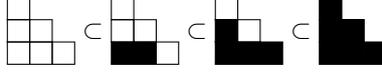
\begin{ytableau}
~  \\
  & \\
&&
  \end{ytableau}\subset
  \ytableausetup
 {smalltableaux,centertableaux}
\begin{ytableau}
~  \\
  & \\
*(black)&*(black)&
  \end{ytableau}\subset
\begin{ytableau}
~  \\
  *(black)& \\
*(black)&*(black)&*(black)
  \end{ytableau}
\subset
\begin{ytableau}
*(black)  \\
  *(black)&*(black) \\
*(black)&*(black)&*(black)
  \end{ytableau}\]
  is \emph{not} a $1$-flag, as condition (iii) is not satisfied by $R_1\subset R_2$.

 \item The flag of regions ($s=3$)\[
\ytableausetup
 {smalltableaux,centertableaux}
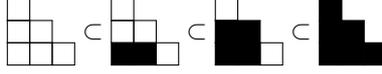
\begin{ytableau}
~  \\
  & \\
&&
  \end{ytableau}\subset
  \ytableausetup
 {smalltableaux,centertableaux}
\begin{ytableau}
~  \\
  & \\
*(black)&*(black)&
  \end{ytableau}\subset
\begin{ytableau}
~  \\
  *(black)&*(black) \\
*(black)&*(black)&
  \end{ytableau}
\subset
\begin{ytableau}
*(black)  \\
  *(black)&*(black) \\
*(black)&*(black)&*(black)
  \end{ytableau}\]
  is \emph{not} a $1$-flag, as condition (ii) is not satisfied by $R_2$ (region not saturated in the south-east direction).
\eite\eexem

 Consider a $1$-flag $R_\bullet$ and the associated $1$-chain $\mathcal E_\bullet$ (again, the pair $(\mathcal F,\theta)$ is fixed for now). We will denote by $|\mathcal E_k|=\# S_k$ the number of boxes in $S_k$ and set $(n_k,p_k)=[\mathcal E_k]$.
 Thanks to Definition~\ref{1def}, each $1$-flag can be seen as the permutation $\sigma$ on $s$ elements given by\[
 (\sigma(1),\dots,\sigma(s))=(|\mathcal E_1|,\dots,|\mathcal E_{s}|).\]
 It will be convenient to represent $1$-flags with horizontal strips, for instance write \[
R_\bullet=\ytableausetup
 {smalltableaux,centertableaux}
\begin{ytableau}
~&&&&&  \\
   \\
&&&\\
& &  &&\\
&  \\
&&
  \end{ytableau}\]
for the $1$-flag associated to $\sigma=(3,2,5,4,1,6)$.

Take two integers $1\le k<j \le s$. Note that $n_j<n_k$ implies $|\mathcal E_j|<|\mathcal E_k|$.
When $n_j<\min\{n_k,\ldots,n_{j-1}\}$, the $\aaa_\text{Higgs}$-slope of the chain (see~\cite[Definition 2.10]{BGGH})\[
\mathcal E_{s}\rarrow\dots\rarrow\mathcal E_j=\dots=\mathcal E_j\rarrow\mathcal E_{k-1}\rarrow\dots\rarrow\mathcal E_1
\]
is the $\boldsymbol\aaa$-slope of the (noncanonical!)\ subregion obtained by only considering the $\# S_j$ leftmost boxes in each $S_t$, $k\le t\le j$. We denote by $R_k^j$ the complementary of this region. For instance if $s=6$ and $\sigma=(3,2,5,4,1,6)$, and if we represent the strips $S_t$ horizontally, we get \[
R_\bullet=\ytableausetup
 {smalltableaux,centertableaux}
\begin{ytableau}
~&&&&&  \\
   \\
&&&\\
& &  &&\\
&  \\
&&
  \end{ytableau}\quad\Rarrow\quad
R_3^4=\ytableausetup
 {smalltableaux,centertableaux}
\begin{ytableau}
~&&&&&  \\
   \\
&&&\\
& &  &&*(black)\\
&  \\
&&
  \end{ytableau},\qquad
R_2^5=\ytableausetup
 {smalltableaux,centertableaux}
\begin{ytableau}
~&&&&&  \\
   \\
&*(black)&*(black)&*(black)\\
& *(black)& *(black) &*(black)&*(black)\\
& *(black) \\
&&  \end{ytableau}
  \]
    where the black boxes are the ones contained in the region.
Note that in order to take  the $\boldsymbol\aaa$-slope, one has to project the boxes in the south direction, and then proceed to the previously mentioned filling~(\ref{TS}) of $T_s$. For instance in our example $\mu^{\boldsymbol\aaa}(R_3^4)=\mu(\aaa_2(-l))$.

  Similarly, when $n_k<\min\{n_{k+1},\dots,n_j\}$, the $\aaa_\text{Higgs}$-slope of the chain\[
\mathcal E_{s}\rarrow\dots\rarrow\mathcal E_{j+1}\rarrow\mathcal E_k=\dots=\mathcal E_k\rarrow \dots\rarrow\mathcal E_1
\]
is the $\boldsymbol\aaa$-slope of the subregion obtained by only considering the $\# S_k$ leftmost boxes in each $S_t$, $k\le t\le j$. We denote by $\check R_k^j$ the complementary of this region, and with the same example $\sigma=(3,2,5,4,1,6)$, we have\[
\check R_2^4=\ytableausetup
 {smalltableaux,centertableaux}
\begin{ytableau}
~&&&&&  \\
   \\
&&*(black)&*(black)\\
& & *(black) &*(black)&*(black)\\
&  \\
&&  \end{ytableau},\qquad
\check R_5^6=\ytableausetup
 {smalltableaux,centertableaux}
\begin{ytableau}
~&*(black)&*(black)&*(black)&*(black)&*(black)  \\
   \\
&&&\\
& &&&\\
& \\
&&  \end{ytableau},
  \]
  where $\mu^{\boldsymbol\aaa}(\check R_5^6)=\sum_{1\le t\le5}\aaa_t$.
  The article~\cite{BGGH} gives necessary and sufficient conditions (C0,C1,C2,C3) for a type $(n_\bullet,p_\bullet)$ to be semistable, meaning that generically, chains of this type are semistable with respect to the $\aaa^\text{Higgs}$-slope. In our context, it yields the following.

 \prop\label{recap}
A $1$-chain $\mathcal E_\bullet$ is of semistable type if it satisfies \[
n_{i-1}=n_i\Rarrow p_i\le p_{i-1}\qquad (\textup C_0)\]
 and if for any $1\le k<j\le s$, the associated $1$-flag $R_\bullet$ satisfies\[
  \mu^{\boldsymbol\aaa}(R_{k})&\le\mu \qquad (\textup C_k)\\
n_j<\min\{n_k,\ldots,n_{j-1}\}\Rarrow \mu^{\boldsymbol\aaa}(R_k^j)&\ge\mu \qquad (\textup C_k^j)\\
 n_k<\min\{n_{k+1},\dots,n_j\}\Rarrow\mu^{\boldsymbol\aaa}(\check R_k^j)&\le\mu \qquad (\check{\textup C}_k^j) \]
 where $\mu=d/r=\mu(\mathcal F)$.
 \eprop

\rema~

\bite
\item Note first that $(\ref{slope})\Rarrow(\textup C_k)$ for any $1$-flag $R_\bullet$.
\item  Also, if $\sigma$ is the permutation associated to a $1$-chain $\mathcal E_\bullet$ of type satisfying  $n_{i-1}=n_i$ but $p_{i-1}<p_i$, one can always multiply on the left $\sigma$ by the transposition $(i-1,i)$ in order to satisfy the condition (C0) - without impacting any of the other conditions.
\eite\erema

 We are going to construct a $1$-chain satisfying this set of conditions, under the assumption~(\ref{slope}).
 Consider a $1$-flag $R_\bullet$, and denote by $\mathcal E_\bullet$ and $\sigma$ the chain and permutation associated to $R_\bullet$, and set\[
 \rho_k^j\sigma=(|\mathcal E_1|,\ldots,\widehat{|\mathcal E_k|},\ldots,|\mathcal E_j|,|\mathcal E_k|,\ldots,|\mathcal E_{s}|)\]
 where $~\widehat{~}~$ means that we remove the underneath entry (this is just the multiplication on the right by the cycle $(k,\ldots,j)$). Denote by $\rho_k^j\mathcal E_\bullet$ and $\rho_k^jR_\bullet$ the corresponding chain and flag. If $\sigma=(3,2,5,4,1,6)$, $k=2$, $j=4$, we have\[
\rho_2^4R_\bullet=\ytableausetup
 {smalltableaux,centertableaux}
\begin{ytableau}
~&&&&&  \\
   \\
   &  \\
&&&\\
& & &&\\
&&  \end{ytableau}.
  \]
The point is that if $(\check{\text C}_k^j)$ fails for $R_\bullet$, then it doesn't for $\rho_k^jR_\bullet$, and if $({\text C}_k^j)$ fails for $R_\bullet$, then it doesn't for $(\rho_k^j)^{-1}R_\bullet$. The following is true thanks to~\ref{applemm}.

\prop\label{algo}
For any $R_\bullet$, there exists a finite sequence $((k_1,j_1,\epsilon_1),\dots,(k_t,j_t,\epsilon_t))$ of integers such that \benu[(i)]
\item $1\le k_p<j_p$ and $\epsilon_p=\pm1$ for all $p=1\dots t$;
\item  if $(\check{\text C}_{k_p}^{j_p})$ (resp. $({\text C}_{k_p}^{j_p})$) fails for $(\rho_{k_{p-1}}^{j_{p-1}})^{\epsilon_{p-1}}\dots(\rho_{k_1}^{j_1})^{\epsilon_1}R_\bullet$ then $\epsilon_{p}=1$ (resp $-1$) for all $p=1\dots t$;
\item $\rho R_\bullet:=(\rho_{k_{t}}^{j_{t}})^{\epsilon_t}\dots(\rho_{k_1}^{j_1})^{\epsilon_1}R_\bullet$ satisfies all conditions $({\text C}_k^j)$ and $(\check{\text C}_k^j)$.
\eenu
\eprop

As explained at the beginning of the section, the following completes the proof of Theorem~\ref{result2}.

\prop\label{1chain} There is an injective map $\boldsymbol \kappa$ from the set of Jordan types $\boldsymbol\aaa\vdash(r,d)$ satisfying~\textup{(\ref{slope})} to the set of semistable types of chains $\mathcal E_\bullet$ such that $\sum_k[\mathcal E_k((k-1)\Omega)]=(r,d)$.\eprop

\demo
Thanks to~\ref{algo}, from any Jordan type $\boldsymbol\aaa\vdash(r,d)$, we build a $1$-flag satisfying all conditions $({\text C}_k^j)$ and $(\check{\text C}_k^j)$. Thanks to~\ref{recap}, the associated $1$-chain is of semistable type.
Denote by $(n_\bullet,p_\bullet)$ this type, the map $\boldsymbol \kappa:\boldsymbol\aaa\mapsto (n_\bullet,p_\bullet)$ hence built is injective because any type of $1$-chain characterizes $\boldsymbol\aaa$, according to the following argument.

We prove more generally that for any $s$ and any $1$-flag $R_\bullet$ of $T_s$, the map ${\boldsymbol\iota}: \boldsymbol\aaa\mapsto[\mathcal E_\bullet]$ from Jordan types of length $s$ to types of chains of length $s$ is injective, where $\mathcal E_k=\mathcal F_{R_k}/\mathcal F_{R_{k-1}}(-(k-1)\Omega)$ for any choice of Higgs pair $(\mathcal F,\theta)$ of type $\boldsymbol\aaa$ (${\boldsymbol\iota}$ doesn't depend on this choice). We use induction on $s$, initialization being trivial at $s=0$. Then take $s>0$, $\boldsymbol\aaa$ of length $s$, and $\mathcal E_\bullet$ as above. Consider $t$ maximizing $\operatorname{rk}(\mathcal E_k)$, as well as $\deg(\mathcal E_k)$ if several pieces of the chain $\mathcal E_\bullet$ have same rank. Then \begin{equation}\label{here}
\alpha_1+\alpha_2(-l)+\dots+\aaa_s(-(s-1)l)=[\mathcal F_{R_t}]-[\mathcal F_{R_{t-1}}](l)=[\ker\theta].\end{equation}
Consider the $1$-flag $\bar R_\bullet$ of $T_{s-1}$ obtained by quotienting by $\ker\theta$, then the associated map $\bar{\boldsymbol\iota}$ is injective by induction hypothesis. Hence, writing $\bar{\boldsymbol\aaa}=(\aaa_2,\dots,\aaa_s)$, we get\[{\boldsymbol\iota}(\boldsymbol\aaa)=\boldsymbol\iota(\boldsymbol\aaa')\Rarrow\bar{\boldsymbol\iota}(\bar{\boldsymbol\aaa})=\bar{\boldsymbol\iota}(\bar{\boldsymbol\aaa}')\Rarrow\bar{\boldsymbol\aaa}=\bar{\boldsymbol\aaa}'\Rarrow{\boldsymbol\aaa}={\boldsymbol\aaa}',\]
 thanks to~\ref{here} for the last implication.

For instance here with $t=s=6$,
\[
R_\bullet=\ytableausetup
 {smalltableaux,centertableaux}
\begin{ytableau}
~&&&&& \times \\
   \\
&&&\\
& &  &\times&\times\\
&  \\
\times&\times&\times
  \end{ytableau}\Rarrow
 \bar R_\bullet=\ytableausetup
 {smalltableaux,centertableaux}
\begin{ytableau}
~&&&&\\
   \\
&&&\\
& & \\
&    \end{ytableau}\]
where $\times$'s correspond to $\ker\theta$.
\edemo

One could wonder how geometric this bijection is. We state the following.

\theo For any semistable Jordan type $\boldsymbol\aaa$, we have $\overline{ \mathbf\Lambda_{\boldsymbol{\aaa}}}=\overline{\mathbf C_{\boldsymbol \kappa(\boldsymbol\aaa)}^-}$.
\etheo

\demo Denote by $\boldsymbol \lambda$ the bijection from the set of semistable Jordan types to the set of semistable types of chains given by\[
\overline{ \mathbf\Lambda_{\boldsymbol{\aaa}}}=\overline{\mathbf C^-_{\boldsymbol \lambda (\boldsymbol\aaa)}}.\]
We want to prove that $\boldsymbol \lambda=\boldsymbol\kappa$. We will use the partial order $\le$ on Jordan cells introduced in~\cite[(2.8)]{SaSchi}:\[
\boldsymbol\bbb\le\boldsymbol\aaa&\Leftrightarrow \forall k : [\ker\theta^k]\le[\ker\theta'^k]\\
& \Leftrightarrow \forall k : [\ker\theta'^k]-[\ker\theta^k]\in\mathbb H.\]

 Consider a semistable type $\boldsymbol\aaa$ and set $\boldsymbol\bbb=\boldsymbol\lambda^{-1}\boldsymbol\kappa(\boldsymbol\aaa)$.
 Consider $(\mathcal F,\theta)\in\mathbf\Lambda_{\boldsymbol\bbb}$ such that there exists a $\theta$-stable flag $(\mathcal F_1\subset\dots\subset\mathcal F_s=\mathcal F)$ of $\mathcal F$ of graded type $\boldsymbol\kappa(\boldsymbol\aaa)$: if $\mathcal E_k=\mathcal F_k/\mathcal F_{k-1}(-(k-1)\Omega)$, we have\[
[\mathcal E_\bullet]=\boldsymbol\kappa(\boldsymbol\aaa).\]
By definition, such a pair $(\mathcal F,\theta)$ exists in $\mathbf\Lambda_{\boldsymbol\bbb}^\text{sst}$, and $\lim_{t\to0}t.(\mathcal F,\theta)=(\mathcal E_\bullet,\phi_\bullet)$ is a semistable chain, where $\phi_\bullet$ is induced by $\theta$.

Assume that $\boldsymbol\bbb\neq\boldsymbol\aaa\vdash(r,d)$, and since we know  that for a fixed pair $(r,d)\in \mathbb H$ there are finitely many semistable Jordan types of size $(r,d)$, pick $\boldsymbol\aaa$ $\le$-maximal of size $(r,d)$ with that property.
 We conclude thanks to the following Lemma~\ref{enfin?}. Indeed by maximality, $\boldsymbol\bbb>\boldsymbol\aaa$ would imply $\boldsymbol\kappa(\boldsymbol\bbb)=\boldsymbol\lambda(\boldsymbol\bbb)=\boldsymbol\kappa(\boldsymbol\aaa)$ hence again $\boldsymbol\aaa=\boldsymbol\bbb$ by injectivity of $\boldsymbol\kappa$.
\edemo

  \lemm\label{enfin?}
We have $\boldsymbol\bbb\ge\boldsymbol\aaa$.
\elemm

\demo
Fix an arbitrary auxiliary pair $(\mathcal G,\psi)$ of type $\boldsymbol\aaa$, and denote by $\mathcal G_\bullet$ the flag constructed along the proof of~\ref{1chain}, whose associated chain is of type $\boldsymbol\kappa(\boldsymbol\aaa)$.
The pair $(\mathcal F,\theta)\in \mathbf\Lambda_{\boldsymbol\bbb}$ being as above, we want to prove that $[\ker(\psi^k)]\le[\ker(\theta^k)]$ for every $k$. Note that $[\mathcal F_k]=[\mathcal G_k]$ for every $k$ and that $[\ker(\psi^k)]$ only depends on $\boldsymbol\aaa$. We prove the following statement, which doesn't involve semistability.

\bigskip

\textit{Consider $(\mathcal F,\theta)\in \mathbf\Lambda_{\boldsymbol\bbb}$ and $(\mathcal G,\psi)\in \mathbf\Lambda_{\boldsymbol\aaa}$, a $1$-flag $\mathcal G_\bullet$ and a flag $\mathcal F_\bullet$ such that $\theta$ induces a chain $(\mathcal E_\bullet,\phi_\bullet)$ of same type as $gr(\mathcal G_\bullet)$. Then $[\ker(\psi^k)]\le[\ker(\theta^k)]$ for every $k$.}

\bigskip

We prove it by induction on the common length $s$ of $\mathcal E_\bullet,gr(\mathcal G_\bullet)$.
There is nothing to prove for $s=0$. Assume that it is known for all lengths $<s$.
Consider for any $k$ the smallest index $t_k$ such that $\ker(\psi^k)\subseteq\mathcal G_{t_k}$.

We prove by induction on $k\ge0$ that $[\ker(\theta^k)]\ge[\ker(\psi^k)]$ together with the following refinement \[
t_{k+1}>t_k&\Rarrow[\ker(\psi^k)]\underset{\eta}{\le}[\mathcal F_{t_{k+1}-1}\cap\ker(\theta^k)]\underset{\eta'}{\le}[\mathcal F_{t_{k+1}}\cap\ker(\theta^k)]\underset{\eta''}{\le}[\ker(\theta^k)]\\
t_{k+1}=t_k&\Rarrow[\ker(\psi^k)]\underset{\eta'}{\le}[\mathcal F_{t_{k+1}}\cap\ker(\theta^k)]\underset{\eta''}{\le}[\ker(\theta^k)]\]
where we mean $[\mathcal B]-[\mathcal A]=\gamma\in\mathbb H$ by $[\mathcal A]\underset{\gamma}{\le}[\mathcal B]$. There is nothing to prove for $k=0$.

Assume that the property is true for some $k$, and that $t_{k+1}>t_k$. We have, using the notation $f_p=[\mathcal F_p]=[\mathcal G_p]$,\[
\left[\dfrac{\ker(\theta^{k+1})}{\ker(\theta^k)}\right]&\ge\left[\dfrac{\mathcal F_{t_{k+1}}\cap\ker(\theta^{k+1})}{\mathcal F_{t_{k+1}}\cap\ker(\theta^k)}\right]\\
&=\left[\ker\left(\dfrac{\mathcal F_{t_{k+1}}}{\mathcal F_{t_{k+1}}\cap\ker(\theta^k)}\to\dfrac{\mathcal F_{t_{k+1}-1}}{\mathcal F_{t_{k+1}-1}\cap\ker(\theta^k)}(\Omega)\right)\right]\\
&\ge f_{t_{k+1}}-([\ker(\psi^k)]+\eta+\eta')-f_{t_{k+1}-1}(l)+([\ker(\psi^k)]+\eta)(l)\\
\Rarrow [\ker(\theta^{k+1})]&\ge\eta''+\eta(l)+f_{t_{k+1}}-(f_{t_{k+1}-1}-[\ker(\psi^k)])(l)
\]
but $t_{k+1}>t_k$ ensures that $f_{t_{k+1}}-(f_{t_{k+1}-1}-[\ker(\psi^k)])(l)=[\ker(\psi^{k+1})]$ so we get $ [\ker(\theta^{k+1})]\ge [\ker(\psi^{k+1})]$.

Assume now that $t_{k+1}=t_k$ and set $\eta=[\ker(\psi^k)]-[\mathcal F_{t_{k+1}-1}\cap\ker(\theta^k)]\in\ZZ^2$. We get this time\[
 [\ker(\theta^{k+1})]&\ge\eta''-\eta(l)+f_{t_{k+1}}-(f_{t_{k+1}-1}-[\ker(\psi^k)])(l).
\]
There exists some $r\ge0$ such that \[
[\ker(\psi^k)]-[\mathcal G_{t_{k+1}-1}\cap\ker(\psi^k)]= \aaa_k+\dots+\aaa_{k-r}
\]
which implies
\[
f_{t_{k+1}}-(f_{t_{k+1}-1}-[\ker(\psi^k)])(l)=[\ker(\psi^{k+1})]+\alpha_k(l)+\dots+\aaa_{k-r}(l).\]
We use our induction hypothesis on the chains induced on $\mathcal F_{t_{k+1}-1},\mathcal G_{t_{k+1}-1}$ to get
\[
\aaa_k+\dots+\aaa_{k-r}-\eta&=\aaa_k+\dots+\aaa_{k-r}-[\ker(\psi^k)]+[\mathcal F_{t_{k+1}-1}\cap\ker(\theta^k)]\\
&\ge \aaa_k+\dots+\aaa_{k-r}-[\ker(\psi^k)]+[\mathcal G_{t_{k+1}-1}\cap\ker(\psi^k)]\\
&=0\]
and thus $[\ker(\theta^{k+1})]\ge [\ker(\psi^{k+1})]$ again.

Assume $t_{k+2}>t_{k+1}$. If $t_{k+1}>t_{k}$, we have from the above chains of inequalities \[
[\ker(\psi^{k+1})]&\le[\mathcal F_{t_{k+1}}\cap\ker(\theta^{k+1})]-\eta(l)\\
&\le [\mathcal F_{t_{k+1}}\cap\ker(\theta^{k+1})]\\
&\le [\mathcal F_{t_{k+2}-1}\cap\ker(\theta^{k+1})].\]
If $t_{k+1}=t_{k}$ we have \[
[\mathcal F_{t_{k+2}-1}\cap\ker(\theta^{k+1})]&\ge [\mathcal F_{t_{k+1}}\cap\ker(\theta^{k+1})]\\
&\ge\ker(\psi^{k+1})-\eta(l)+\alpha_k(l)+\dots+\aaa_{k-r}(l)\\
&\ge\ker(\psi^{k+1}).\]
Now assume $t_{k+2}=t_{k+1}$. We have $[\ker(\psi^{k+1})]\le [\mathcal F_{t_{k+1}}\cap\ker(\theta^{k+1})]$ for the same reasons as above, which concludes the proof.
\edemo

\section{Polytopal description}

We first reformulate Theorem~\ref{result2}, in order to use a slightly different set of inequalities to characterize the semistability. For any region $R\in\mathcal R$ we denote by $R_k$ the height of its $k$-th column, counted from the right. For instance \[
R=\ytableausetup
 {smalltableaux,centertableaux}
\begin{ytableau}
  \\
 &  \\
&&\\
*(black)& *(black)&  &\\
*(black)&*(black)&*(black)&  *(black) &
  \end{ytableau}\Rightarrow (R_5,R_4,R_3,R_2,R_1)=(2,2,1,1,0).\]
We set $\mathcal R_p=\{R\in\mathcal R\mid R_1=\ldots=R_{p-1}=0\neq R_p\}$ and $\mathcal R_{>1}=\cup_{p>1}\mathcal R_p$. Note that from the definition of $\mathcal R$, there is a bijective map $\mathcal R\rarrow\mathcal R$ sending $R=(R_k)$ to $\bar R=(k-R_k)$.  For instance\[
R=\ytableausetup
 {smalltableaux,centertableaux}
\begin{ytableau}
  \\
 &  \\
&&\\
*(black)& *(black)&  &\\
*(black)&*(black)&*(black)&  *(black) &
  \end{ytableau}\mapsto
\bar  R=\ytableausetup
 {smalltableaux,centertableaux}
\begin{ytableau}
  \\
 &  \\
*(black)&&\\
*(black)& *(black)& *(black) &\\
*(black)&*(black)&*(black)&  *(black) &  *(black)
  \end{ytableau}\]
and we see that this map swaps $\mathcal R_1$ and $\mathcal R_{>1}$. For any Jordan type $(\mathbf r,\mathbf d)\vdash(r,d)$ and any region $R\in\mathcal R$, we finally define\[
p_R(\mathbf r)=\dfrac{ r_R}{r}=\dfrac{\sum_{\blacksquare\in R}\operatorname{rank}\blacksquare}{r}=\dfrac{\sum_kR_kr_k}{r}.
\]

\prop\label{ineq} A Jordan type $\boldsymbol\aaa=(\mathbf r,\mathbf d)\vdash(r,d)$ is semistable if and only if for any region $R\in\mathcal R_{>1}$ we have \begin{align}\label{slope2}
p_R(\mathbf r)d+l\sum_k\frac{R_k(R_k-1)}{2}r_k\le\sum_kR_kd_k\le p_R(\mathbf r)d+l\sum_k\frac{R_k(\bar R_k+k-1)}{2}r_k.\end{align}
\eprop
 \rema Note that the bounds do not depend on $\mathbf d$ and that $d_1$ never appears in the central term. Hence, because of~(\ref{deun}), it is only subject to\begin{align}\label{remd1}
d_1=d-\sum_{k\ge2}kd_k+l\sum_k\frac{k(k-1)}{2}r_k.\end{align}
\erema
\demo First note that \[
 d_R=\sum_{\blacksquare\in R}\deg\blacksquare=\sum_kR_kd_k-l\sum_k\frac{R_k(\bar R_k+k-1)}{2}r_k\]
so that~(\ref{slope}) is equivalent to\[
\sum_kR_kd_k\le p_R(\mathbf r)d+l\sum_k\frac{R_k(\bar R_k+k-1)}{2}r_k.\]
Also,\[
 d_{\bar R}=d- d_R-l\sum_kR_k\bar R_kr_k\]
so that~(\ref{slope}) with respect to $\bar R$ is equivalent to\[
&d- d_R-l\sum_kR_k\bar R_kr_k\le p_{\bar R}(\mathbf r)d=(1-p_R(\mathbf r))d\\
\Leftrightarrow\quad&p_R(\mathbf r)d+l\sum_k\frac{R_k(R_k-1)}{2}r_k\le \sum_kR_kd_k.\]
This concludes the proof since $\mathcal R=\mathcal R_1\sqcup\mathcal R_{>1}$.
\edemo

Assume that $r>0$. Then if $(\mathcal F,\theta)\in\boldsymbol\Lambda_{r,d}^\textup{sst}$ has type $\boldsymbol\aaa$ of length $s$, we necessarily have $r_s>0$ since $r_s=\operatorname{rank}\mathcal F_{s-1}$ and $\mathcal F_{s-1}\subset\mathcal F$.

\coro\label{polytop} Fix a partition $\mathbf r\vdash r$ of length $s$. The set of degree vectors $\mathbf d=(d_1,\ldots,d_s)$ such that $(\mathbf r,\mathbf d)\vdash(r,d)$ is semistable is the intersection $\mathcal P_{\mathbf r,d}^\ZZ$ of the integral lattice $\ZZ^s$ with a convex $(s-1)$-polytope $\mathcal P_{\mathbf r,d}$.\ecoro

\demo The polytope $\mathcal P_{\mathbf r,d}$ is defined by the set of linear inequalities~(\ref{slope2}) together with the hypothetical extra conditions $d_k\ge0$ every time we have $r_k=0$. It is $(s-1)$-dimensional because of~(\ref{remd1}). Also because of this equation, note that a facet may be given by the equation\begin{align}\label{extrafacet}
\sum_{k\ge2}kd_k= d+l\sum_k\frac{k(k-1)}{2}r_k\end{align}
if $r_1=0$.\edemo

We get the following correspondence between polytopes associated to different degrees.
\prop\label{indep} Consider $d,d'\in\ZZ$ and fix a partition $\mathbf r\vdash r$ of length $s$. The translation $\boldsymbol\tau=(d'-d)\mathbf r/r$ of $\RR^s$ induces a bijection $\mathcal P_{\mathbf r,d}\isom\mathcal P_{\mathbf r,d'}$.
\eprop

\demo First rewrite~(\ref{slope2}) in the following way\begin{align}\label{slope3}
b^-_R(\mathbf r)\le \left(\mathbf d-\dfrac{d}{r}\mathbf r\right)\cdot R\le b^+_R(\mathbf r)\end{align}
where\[
b^-_R(\mathbf r)=l\sum_k\frac{R_k(R_k-1)}{2}r_k,\quad&\quad
b^+_R(\mathbf r)=l\sum_k\frac{R_k(\bar R_k+k-1)}{2}r_k,\\
\mathbf r\cdot R=\sum_k R_kr_k,\quad&\quad
\mathbf d\cdot R=\sum_k R_kd_k.\]
Note that since we consider $R\in\mathcal R_{>1}$ we can replace $\mathbf d$ by $\mathbf d^*=(0,d_2,\ldots,d_s)$ everywhere.

Finally, also note that if $r_1=0$ the equation~(\ref{extrafacet}) fits in this study since it can be written\[
\left(\mathbf d^*-\frac{d}{r}\mathbf r\right)\cdot T_s=l\sum_k\frac{k(k-1)}{2}r_k\]
where $T_s=(k)$ corresponds to the full tableau.\edemo

It is not clear how this bijection restricts to integral points, but we know, as explained in the introduction, that thanks to~\cite{Mellit} we have the following.

\theo The quantity $\sum_{\mathbf r\vdash r}\#\mathcal P^\ZZ_{\mathbf r,d}$ does not depend on degrees $d$ coprime with $r$, and is equal to $A_{g,r}(1)$.\etheo

The Propostion~\ref{indep} together with the rather explicit description~(\ref{slope3}) however shed an interesting light on Mellit's independence result, given our direct geometric approach of the global nilpotent cone, and its combinatorial flavor.

In fact, Rodriguez Villegas defined in~\cite{RVA} a refinement $A_{g,\mathbf r}$ of $A_{g,r}$ for each partition $\mathbf r$ of $r$, satisfying $A_{g,r}(q)=\sum_{\mathbf r\vdash r}A_{g,\mathbf r}(q)$, and establishes closed formulas for the quantities $A_{g,\mathbf r}(1)$. Based on computations for small values of $l(\mathbf r)$ and the case $\mathbf r=(1^r)$ established by Reineke~\cite[\S7]{MR2889742}, the following is expected.

\begin{conj} We have $\#\mathcal P^\ZZ_{\mathbf r,d}=A_{g,\mathbf r}(1)$.
\end{conj}

\appendix

\section{Descent argument, by Anton Mellit}
In this appendix we show that the procedure used to construct a flag satisfying conditions of Proposition \ref{recap} eventually stops. We briefly recall the setup. Fix $l\geq 0$ and suppose $r_1,\ldots,r_s,d_1,\ldots,d_s$ is a semistable Jordan type. Recall that $1$-flags $R_\bullet$ correspond to permutations $(\sigma(1),\ldots,\sigma(s))$ where $\sigma(k)=|S_k|$,  $S_k=R_k\setminus R_{k-1}$. A $1$-flag is visually represented by a collection of horizontal strips of lengths $\sigma(1),\cdots,\sigma(s)$. The $i$-th box in $S_k$ has rank $r_i$ and degree
\[
d_i + l r_i \#\{j<i\,:\, \sigma(j)\geq i\}.
\]
For every subset of boxes $Q$ the rank $r_Q$ and degree $d_Q$ are defined by summing the ranks resp. degrees of its boxes. The slope is the ratio $\mu(Q)=\frac{d_Q}{r_Q}$, which by definition equals $+\infty$ in the case $r_Q=0$. The slope of the Jordan type is the slope of the set of all boxes and denoted by $\mu$. This does not depend on the flag. Semistability means that $\mu(R_k)\leq \mu$ for any $1\leq k\leq s$ and any flag $R_\bullet$.

In the procedure described in Section 3 we consider regions $R_k^j$ resp. $\check R_k^j$ where $1\leq k<j\leq s$ and $n_j<\min\{n_{k},\ldots,n_{j-1}\}$ resp. $n_k<\min\{n_{k+1},\ldots,n_j\}$ and check the conditions ${\textup C}_k^j$ resp. $\check{\textup C}_k^j$, see Proposition \ref{recap}. If $\check{\textup C}_k^j$ fails, we apply permutation $\rho_k^j$, which places $k$-th strip at position $j$ and shifts all the strips in between down. The resulting flag is denoted $\rho_k^j R_\bullet$. If on the other hand ${\textup C}_k^j$ fails, we apply the inverse permutation by placing the $j$-th strip at position $k$ and shifting all the strips in between up. The resulting flag is denoted $(\rho_k^j)^{-1} R_\bullet$.

\begin{lemme}\label{applemm}
    In the above setup, there is a function $f$ from the set of $1$-flags to the set of real numbers such that if the condition $(\check{\textup C}_k^j)$ fails for some $k$ and $j$, then $f(\rho_k^j R_\bullet)< f(R_\bullet)$, and if the condition $({\textup C}_k^j)$ fails for some $k$ and $j$, then $f((\rho_k^j)^{-1} R_\bullet)< f(R_\bullet)$. Therefore the procedure eventually stops.
\end{lemme}
\begin{proof}
    Suppose $l>0$ first. To any $1$-flag $R_\bullet$ we associate the following real number:
\begin{equation}\label{the function}
f(R_\bullet) := \sum_{i=1}^s  (d_{S_i} - \mu\, r_{S_i})^2.
\end{equation}
 Set $\delta=  f(\rho_k^j R_\bullet) - f(R_\bullet) $. We have
        \begin{align*}
    \delta &=(d_{S_k} + (j-k) l r_{S_k} - \mu\, r_{S_k})^2 + \sum_{i=k+1}^j (d_{S_i} - l r_{S_k} - \mu\, r_{S_i})^2
    \\
&\qquad\qquad\  -\sum_{i=k}^j  (d_{S_i} - \mu\, r_{S_i})^2  \\
    &= 2 (j-k) l r_{S_k}( d_{S_k}- \mu\, r_{S_k}) + \left((j-k) l r_{S_k}\right)^2 - \sum_{i=k+1}^j\left( 2 l r_{S_k} (d_{S_i} - \mu\, r_{S_i}) +  l^2 r_{S_k}^2\right)\\
    &= 2 l r_{S_k} \left((j-k) ( d_{S_k}- \mu\, r_{S_k}) + \binom{j-k+1}2 l r_{S_k} - \sum_{i=k+1}^j(d_{S_i} - \mu\, r_{S_i}) \right).
        \end{align*}
    
    Note that $r_{S_k}>0$, otherwise by placing $S_k$ at the first position we would obtain a flag with $\mu(R_1)=+\infty$ which contradicts semistability. So the overall sign coincides with the sign of the expression in parenthesis. The latter expression can be simplified as follows:
    \[
    \mu \sum_{i=k+1}^j (r_{S_i} - r_{S_k}) - \sum_{i=k+1}^j\left(d_{S_i} - d_{S_k} - (i-k) l r_{S_k}\right) = \mu r_{\check R_k^j} - d_{\check R_k^j}.
    \]
    So we have established
    \begin{equation}\label{cool formula}
    f(\rho_k^j R_\bullet) - f(R_\bullet) = 2 l r_{S_k} (\mu r_{\check R_k^j} - d_{\check R_k^j}).
    \end{equation}
    If condition $(\check{\textup C}_k^j)$ fails for $R_\bullet$, then we have $\mu^{\boldsymbol\aaa}(\check R_k^j)>\mu$ and so $f(\rho_k^j R_\bullet) - f(R_\bullet)<0$. If on the other hand condition $({\textup C}_k^j)$ fails for $R'_\bullet=\rho_k^j R_\bullet$, then we have $\mu(R'^j_k) = \mu(\check R_k^j)<\mu$, which implies $f(R'_\bullet)-f((\rho_k^j)^{-1}R'_\bullet)=f(\rho_k^j R_\bullet) - f(R_\bullet)>0$.

    Strictly speaking, we do not need the case $l=0$ because it corresponds to $g=1$, but we include it here for completeness. We construct $f$ using a trick. Let us introduce $l$ as a parameter, and make all $d_Q$ of all regions $Q$ a function of $l$. Then \eqref{the function} produces a family of values $f_l(R_\bullet)$ depending on $l$ for any $1$-flag $R_\bullet$. Let
    \[
    f(R_\bullet) = \left(\frac{\partial}{\partial l} f_l(R_\bullet)\right)\Bigg|_{l=0}.
    \]
    Differentiating the final formula \eqref{cool formula} by $l$ we obtain
    \[
    f(\rho_k^j R_\bullet) - f(R_\bullet) = 2 r_{S_k} (\mu r_{\check R_k^j} - d_{\check R_k^j}),
    \]
    and the proof continues in the same way as for $l>0$.
\end{proof}

\thanks{\textsc{
\noindent
\\
Institut Camille Jordan - Universit\'e Lyon 1}
 \\\small{B\^at. Braconnier,
43 Boulevard du 11 Novembre 1918, 69622 Villeurbanne cedex, France.
}
\\ \textit{e-mail}:\;\texttt{bozec@math.univ-lyon1.fr}}

\end{document}